\documentclass[11pt]{amsart}
\usepackage{amsmath,amsthm,epsfig,amssymb}
\usepackage{amsthm,amsfonts,amssymb,amscd}
\usepackage[utf8x]{inputenc}
\usepackage{textcomp}
\usepackage{tikz}
\usetikzlibrary{positioning, matrix, arrows}
\usepackage[all]{xy}
\usepackage[colorlinks=true,citecolor=blue,urlcolor=blue,linkcolor=blue]{hyperref}
\usepackage{url}

\newtheorem{Theorem}[subsection]{Theorem}
\newtheorem{Proposition}[subsection]{Proposition}

\newtheorem{Lemma}[subsection]{Lemma}

\theoremstyle{definition}
\newtheorem{Remark}[subsection]{Remark}

\newtheorem{proposition-definition}[subsection]{Proposition-Definition}

\title{An  elementary  proof of Lelli Chiesa's theorem on constancy of second coordinate of gonality sequence   }
\author[S.Pal]{Sarbeswar Pal}
\address{IISER - Thiruvananthapuram, Computer Science Building,
College of Engineering Trivandrum Campus,
Trivandrum - 695016, Kerala, India}
\email{spal@iisertvm.ac.in}


\begin{document}

\begin{abstract}
 Let $X$ be a K3 surface and $L$ be an ample line bundle on it. In this article we will give an alternative and elementary proof of Lelli Chiesa's Theorem in the case of $r= 2$. More precisely we will prove that under certain condition the second co-ordinate of the gonality sequence is constant along the smooth curves in the linear system $|L|$.  Using Lelli Chiesa's theorem for $r \ge 3$ we also  extend Lelli Chiesa's Theorem in the case of $r= 2$ in  weaker condition. 
\end{abstract}

\maketitle

\section{Introduction}
Given a smooth irreducible projective curve $C$ and an integer $r$ one can associate an integer $d_r$ as the minimal degree of a line bundle  with $r+1$ sections. Thus to each curve one can associate a sequence $(d_1, d_2, ...)$ called gonality sequence. The first co-ordinate of the gonality sequence is known as the gonality of $C$.
Let $X$ be a smooth projective K3 surface over the field of complex numbers and $L$ be a line bundle on $X$. Then the natural question one can ask whether the gonality sequence remains constant as $C$ varies in $|L|_s$, where $|L|_s = \{C \in |L| : C \text{ is smooth }\}$. The answer of the question is negative. In fact,  Donagi and Morrison pointed out the following easy counter example showing that even the first co-ordinate is not constant.

\textbf{Example}(\cite{3}, 2.2). Let $\pi: X \to \mathbb{P}^2$ be a K3 surface obtained as a double cover of $\mathbb{P}^2$ ramified at a smooth sextic curve. Let $L = \pi^*(\mathcal{O}_{\mathbb{P}^2}(3))$. The general curve of $|L|$ is a plane sextic and hence they have gonality $5$. On the other hand, 
$|L|$ contains a subspace of co-dimension $1$ consisting of bielliptic curves which has gonality $4$.

However, Ciliberto and Pareschi proved that if $L$ is an ample line bundle on a K3 surface $X$, such that $X$ and $L$ are not  simultaneously as in the  Donagi-Morrison's example, then gonality remains constant along $|L|_s$ \cite{2}. 

Naturally one could ask about the behavior of the second co-ordinate.
Note that in the Donagi-Morrison's example, the second co-ordinate (which we will call planarity of $C$ and denote by $\mathcal{P}(C)$) is constant.  

Recently Lelli Chiesa  \cite{7} proved that if $C$ is an an ample curve in $X$ with some extra hypothesis and admits a complete $g^r_d$ computing the Clifford index of $C$  then every curve in the linear system $|\mathcal{O}_X(C)|$ admits a complete $g^r_d$. However if the Clifford index is bigger than $2$, then  the extra hypothesis satisfied automatically. Thus if the Clifford index of $C$ is bigger than $2$ and $C$ admits a complete $g^2_d$ computing the Clifford index of $C$, then 
$\mathcal{P}(C)$ is constant as $C$ varies in the linear system $|\mathcal{O}_X(C)|$.  The question of constancy of $\mathcal{P}(C)$ still remains open if $C$ does not admit a complete $g^2_d$ computing the Clifford index of $C$. For example see section \ref{S2}.
 
 In this article we will  give an independent proof for constancy of $\mathcal{P}(C)$ when $C$ admits a complete $g^2_d$ computing the Clifford index and also few cases when it does not admit a complete $g^2_d$ computing the Clifford index. We prove the following Theorem:
\begin{Theorem}\label{TH}
Let $X$ be a smooth projective $K3$ surface over the field of complex numbers. Let $L$ be an ample line bundle on $X$ such that there is no bi-elliptic curve in $|L|$.  Then every smooth curve in the linear system $|L|$ carry a $g^2_{d^\prime}$, if  one of the following holds:\\
(i)there exist an irreducible smooth curve $C$ in the linear system $|L|$ with a complete $g^r_{d^\prime},$ for some $r,  2 \le r \le 3$  
which computes the Clifford index of $C$.\\
(ii) $L^2 \ge 16$ and there exists a smooth curve $C \in |L|$ with  a complete $g^4_{d^\prime},$ which computes the Clifford index of $C$.

In other words, the second co-ordinate of the gonality sequence  of  smooth curves is constant along the linear system $|L|$. 
\end{Theorem}

{\it Notation:} We work throughout over the field $\mathbb{C}$ of complex numbers. 
If $X$ is a smooth, projective variety, we denote by $K_X$ the canonical bundle on $X$. 
For a coherent sheaf $\mathcal{F}$ on $X$, we denote by $H^i(\mathcal{F})$ the 
$i$-th cohomology group of $\mathcal{F}$ and by $h^i(\mathcal{F})$ its (complex) 
dimension. If $V$ is a vector bundle on $X$, we denote by $V^*$ the dual of $V$. For a
sub-scheme $Z \subset X$, we denote by $\mathcal{I}_Z$ the ideal sheaf of $Z$. A line bundle of degree $d$ is called a complete $g^r_d$ on a smooth projective curve $C$ if it has exactly $r+1$ sections. We  denote by $W^r_d(C)$, the  subvariety  of 
$\text{Pic}^d(C)$  whose  support is the set:
\[
 \text{Supp}(W^r_d(C)) = \{L \in \text{Pic}^d(C): h^0(C, L) \ge r+1\}.
\]

If $r=0$ we denote $W^0_d(C)$ simply by $W_d(C)$.
\section{preliminaries}\label{S1}
In this section we recall the basic properties of the bundle $E_{C,  A}$ of Lazarsfeld  \cite{5} and Tyurin  \cite{9}, associated to an irreducible smooth curve $C$ in $X$ and a globally generated  line bundle $A$ and 
the basic definitions of Clifford index and Clifford dimension. 

Let $X$ be a smooth projective $K3$ surface over the field of complex numbers. Let $C$ be an irreducible smooth curve in $X$ and $A$ be a globally generated line bundle on $C$. Viewing $A$ as a sheaf on $X$, consider the evaluation map
\[
H^0(C, A) \otimes \mathcal{O}_X \to A.
\]
Let $F_{C, A}$ be its kernel and $E_{C, A} := F_{C, A}^*$.  Then $F_{C, A}$ fits in the following exact sequence on $X$.
\[
0 \to F_{C, A} \to H^0(C, A) \otimes \mathcal{O}_X \to A \to 0.
\]
It is easy to check that $F_{C, A}$ is locally free. Dualizing the above exact sequence one gets 
\[
0 \to {H^0(C, A)}^* \otimes \mathcal{O}_X \to E_{C, A} \to \mathcal{O}_C(C) \otimes A^* \to 0.
\]
Then it is easy to check the following properties:
\begin{Lemma}\label{L1}
1. Rank of $E_{C, A} = h^0(C, A)$.\\
2. $\text{det}(E_{C, A}) = \mathcal{O}_X(C)$.\\
3. $c_2(E_{C, A}) = \text{deg}(A)$.\\
4. $h^0(X, {E_{C, A}}^*)= h^1(X, {E_{C, A}}^*) = 0$.\\
5. $E_{C, A}$ is generated by its global sections off a finite set.
\end{Lemma}

\subsection{Clifford index}
Let $C$ be a smooth irreducible complex projective curve of genus $g \ge 2$.
Recall that the Clifford index of a line bundle $A$ on $C$ is the integer
\[
\text{Cliff}(A) = \text{deg}(A) - 2r(A),
\]
where $r(A) =h^0(A)- 1$. The Clifford index of $C$ itself is defined to be
\[
\text{Cliff}(C) = \text{min} \{\text{Cliff}(A) | h^0(A) \ge  2, h^1(A)  \ge 2\}.
\]
We say that a line
bundle $A$ on $C$ contributes to the Clifford index of $C$ if $A$ satisfies the
inequalities in the definition of $\text{Cliff}(C)$; it computes the Clifford index of $C$ if in
addition $\text{Cliff}(C) = \text{Cliff}(A)$. 
\begin{Theorem}( M. Green, R. Lazarsfeld \cite{4})
Let $X$ be a complex projective K3 surface, and let $C \subset X $be a smooth
irreducible curve of genus $g \ge 2$. Then
\[
\text{Cliff}(C^{\prime}) =\text{Cliff}(C)
\]
for every smooth curve $C^{\prime} \in |C|$. Furthermore, if $\text{Cliff}(C)$  is strictly less than the
generic value $[\frac{(g-1)}{2}]$, then there is a line bundle $L$ on $X$ whose restriction to
any smooth $C^{\prime} \in |C|$ computes the Clifford index of $C^{\prime}$.
\end{Theorem}
Given  a curve $C$, we define  its Clifford dimension as 
\[
r= \text{min}\{h^0(A) -1 | A \text{ computes the Clifford index of } C\}.
\]
\begin{Proposition}(Ciliberto, Pareschi \cite{2})\label{PA}
Let $C$ be a smooth and irreducible curve of genus $g$ sitting on a K3
surface $X$ as an ample divisor. Then either $C$ is isomorphic to a smooth plane sextic and
$X, \mathcal{O}_X(C)$  are as in Donagi-Morrison's example or the Clifford dimension of $C$ is $1$.
\end{Proposition}

\section{An example }\label{S2}
In this section we will give an example of a curve $C$ in a K3 surface $X$ such that the Clifford index of $C$ is not computed by a $g^2_d$ but a $g^3_d$.  Therefore  we can not use Lelli Chiesa's Theorem to conclude the constancy of the second co-ordinate of the gonality sequence.  However we will see that  the second co-ordinate remains constant along $\mathcal{O}_X(C)$, which gives an example in support of our Theorem \ref{TH}. 

{\bf Example}: Let $X$ be the K3 surface given by a smooth quartic hypersurface in $\mathbb{P}^3$. Let $C$ be a quadric hypersurface section. In other words, $C$ is a complete intersection of two hypersufaces of degree $4$ and $2$ respectively.  Clearly $C$ is an ample curve in $X$. Then we have following facts \cite[p.199, F-2]{1}:\\
$\bullet $ $  W^1_3(C) = \varnothing$\\
$\bullet$ $W^1_4(C) \ne \varnothing$ \\
$\bullet$ $W^3_8(C) \ne \varnothing$\\
$\bullet$ $W^3_8(C)-W_2(C)) \subset W^1_6(C) $\\
$\bullet$ $W^2_7(C) = W^3_8(C) - W_1(C)$.\\
Thus the Clifford index of $C$ is $2$. Since $W^2_7(C) = W^3_8(C)-W_1(C)$ and $W^3_8(C) -W_2(C) \subset W^1_6(C)$, we  have $W^2_6 = \varnothing$. Therefore the Clifford index of $C$ can not be computed by a $g^2_d$. On the other hand, since $W^3_8(C)$ is non-empty, the Clifford index is computed by a $g^3_d$.  It is clear that $\mathcal{P}(C) = 7$ for all smooth curve $C \in |\mathcal{O}_X(C)|$.

 \section{Structure of $E_{C, A}$}
 Let $C$ be a smooth irreducible curve in a K3 surface $X$ and $A$ be a line bundle of minimal degree $d$ with $3$ sections. Clearly such a line bundle is globally generated.
Let $E_{C, A}$ be the vector bundle constructed as in Section \ref{S1}. Then by Lemma \ref{L1} we have, \\
\begin{align*}
\text{rk}(E_{C, A}) = 3, \text{det}(E_{C, A}) = \mathcal{O}_X(C),  c_2(E_{C, A})\\ 
 = d, h^0(X, {E_{C, A}}^*)= h^1(X, {E_{C, A}}^*) = 0
\end{align*}
and $E_{C, A}$ is globally generated off a finite set.\\
The following Proposition is a slight modification of  a result of Donagi-Morrison.

\begin{Proposition}\label{SP1}
 $E_{C, A}$ is not a simple vector bundle, then we have the following possibilities:\\
(1) There exist a base point free line bundle $N$ and a  rank $2$ vector bundle $F$, globally generated off a finite set such that $E_{C, A} = F \oplus N$.\\
(2) There exist a base point free line bundle $N$, a rank $2$ vector bundle $F$ and a finite set $Z \subset X$ such that $E_{C, A}$ sits in the following exact sequence,
\[
0 \to F \to E_{C, A} \to  N \otimes \mathcal{I}_Z \to 0
\]
and we have $h^0(F) \ge h^0(N) \ge 2$.
\end{Proposition}
\begin{proof}
If $E_{C, A}$ is not simple, then there is an endomorphism $\varphi: E_{C, A} \to E_{C, A}$ which is not of the form $c.Id$ for some scalar $c$, where $Id$ denotes the identity morphism. Let $x \in X$ be a point. Consider an eigen value $c$ of the linear map $\varphi_x: {(E_{C, A})}_x \to {(E_{C, A})}_x$. Then the morphism $\psi := \varphi - c.Id$ is a nonzero morphism, which drops rank everywhere. 

Let $F: = \text{Ker}(\psi), N^{\prime} = Im(\psi)$.  If $E_{C, A}$ is decomposable then we are in situation (1). Let us assume $E_{C, A}$ is indecomposable. If  the rank of the endomorphism $\psi$ is $2$, then one can easily see that the rank of $\psi^2$ is $1$.  Thus with out loss of generality we can assume that $\text{rk}(F) = 2$ and we have a
short exact sequence of the form,
\[
0 \to F \to E_{C, A} \to N^{\prime} \to 0.
\]
 
Since $X$ is a surface, any reflexive sheaf over $X$ is  locally free. Thus $F$ is locally free.  \\ 
Note that  $N := {N^{\prime}}^{**}$ is a line bundle and $N^{\prime} = N \otimes \mathcal{I}_Z$, for some finite set $Z \subset X$.
Thus we have a sequence
\[
0 \to F \to E_{C, A} \to N \otimes \mathcal{I}_Z \to 0.
\]
Since $E_{C, A}$ is globally generated off a finite set, $N$ is also globally generated off a finite set. 
 Since a line bundle on a K3 surface  has  no base points outside its fixed component  [Corollary 3.2, \cite{8}],  it is globally generated. Moreover, since $h^0(X, {(E_{C,A})}^*) = 0$, $N$ is non-trivial.  Thus $h^0(N) \ge 2$.  If $\psi^2 \ne 0$ then the sequence splits and again we are in the situation (1).  If $\psi^2 = 0$, then $h^0(F \otimes N^*) > 0$. Therefore, we have $h^0(F) \ge h^0(N) \ge 2$.
\end{proof}
\begin{Remark}\label{RZ}
Note that if we are in the second case, then $h^0(N^* \otimes F) \ne 0$. Thus $F$ and hence $E_{C, A}$ contains a line subbundle $M$ which admits at least $2$ sections. Thus $E_{C, A}$ fit in the following exact sequence,
\[
0 \to M \to E_{C, A} \to F \to 0,
\]
where $F$ is a torsion free sheaf of rank two generated by its global sections off a finite set and we have the following exact sequence
\[
0 \to F \to F^{**} \to S \to 0,
\]
where $F^{**}$ is the double dual of $F$, $S$ is a coherent sheaf of finite length , in particular supported on a zero-dimensional subscheme $Z$.
 Also note that $c_2(F)= c_2(F^{**}) + |Z|$, where $|Z|$ denotes the length of $Z$.
\end{Remark}
 \begin{Lemma}\label{L2}
 If $E$ is a globally generated vector bundle off a finite set and $c_1(E)^2 >0$,  then $c_2(E) \ge 0$.
 \end{Lemma}
 \begin{proof}
 If $E$ is a globally generated vector bundle off a finite set then for a general subspace $V \subset H^0(E)$ of dimension $\text{rk}(E)$, we have the following exact sequence \cite[ See P.18 ]{2}
 \[
 0 \to V \otimes \mathcal{O}_X \to E \to B \to 0
 \]
 where $B$ is a line bundle on a smooth curve $C \subset X$. Dualizing the exact sequence we have,
 \[
 0 \to E^* \to V^* \otimes \mathcal{O}_X \to A \to 0
 \]
 where $A = K_C \otimes B^*$.  If $\text{deg}(A) < 0$, then degree of $B \ge 2g^{\prime}-1$, where $g^{\prime}$ is the genus of $C$ and hence $h^0(B) \ge g^{\prime}$. Thus $h^0(E) \ge g^{\prime}+2$. On the other hand since $c_1(E)^2 > 0$, we have $h^0(E) \le h^0(c_1(E)) = g^{\prime} +1$ [Proposition 1.5, \cite{4}]. Thus $c_2(E) = \text{deg}(A) \ge 0$.
 \end{proof}
 
 \section{Trigonal curve in K3 surface}
 In this section we will prove an interesting property of a trigonal curve in a K3 surface .
 \begin{Theorem}\label{TTA}
 Let $C$ be a trigonal curve of genus $g \ge 5$ in a K3 surface $X$. Then the following holds:
 There exists an irreducible curve $\Delta$ such that ${p}_a(\Delta) = 1$ and $\Delta . C = 3$.
 \end{Theorem}
 \begin{proof}
 Since $C$ is a trigonal curve, its Clifford index is $1$. Note that if a $g^r_d$ computes the Clifford index of $C$, then  $d = 2r+1$. 
 On the other hand, for a trigonal curve  if $d_r$ is the minimal degree of a  line bundle with at least $r+1$ sections, then we have (see \cite[Remark 4.5(b)]{6})
 \begin{equation}\label{EEE}
\begin{split}
 d_r  & =3r , 1 \le r \le [\frac{g-1}{3}], \\
         & = r + g -1 - [\frac{g-r -1}{2}], [\frac{g-1}{3}] < r \le g-1 \\
         & =r + g, r \ge g
  \end{split}
\end{equation}
 Therefore, if a line bundle of degree $2r+1$ has at least $r+1$ sections, then $2r+1= d_r$. Now from the expression of $d_r$ in \ref{EEE}, one can conclude that the possibilities are  $r =1$ and $r= g-2$. In other words,    the Clifford index of a trigonal curve can be computed only by a pencil $L$ and $K \otimes L^*$. \\  
 On the other hand,  there exist a line bundle $M$ on $X$ such that $M{\mid_C}$ computes the Clifford index \cite{4}. Therefore, $h^0(C, M{\mid_C}) =2$  or $h^0(\mathcal{O}_X(C) \otimes M^*){\mid_C}) =2 $.  With out loss of generality we assume $h^0(C, M{\mid_C}) =2$. Since $K_C \otimes (M{\mid_C})^*$ also computes the Clifford index, once can see that $h^0(M \otimes \mathcal{O}(-C)) =0$ and hence $h^0(X, M)= 2$ and  $\text{deg}(M{\mid_C})= M.C = 3$.  Therefore a general curve $\Delta$ in $|M|$ is irreducible and has arithmetic genus $1$.  Also we have $\Delta. C =3$, which conclude the Theorem.
 \end{proof}
 
 \section{Main theorem}
In this section we prove the main theorem. \\
If $X$ and $L$ are as in the Donagi-Morrison's example, then we have seen that the planarity remains constant along the smooth curves in $|L|$. 
Let assume $X$ and $L$ are not as in the Donagi-Morrison's example. Let $C$ be an irreducible smooth curve $C$ in the linear system $|L|$ with a complete $g^r_{d^\prime}$, where $2 \le r \le 4$,  which computes the Clifford index of $C$. It is knowm that the gonality is constant along the smooth curves in the linear system $|L|$ \cite{2}. Let $d$ be the gonality.  Also we have the Clifford dimension of every curve in the linear system $|L|$ is $1$ \cite{2}.  Thus the Clifford index of every curve is $d-2$.

 \begin{proof}{ of main Theorem }:\\
 {\bf Case I: r=2}

 Let $C \in |L|$ be a smooth curve. If $C$ is hyperelliptic then the Theorem holds trivially. We assume $C$ is not hyperelliptic.  Let $A$ be  a complete $g^2_{d^\prime}$,  on $C$, computing the Clifford index.  Therefore,  the degree $d^\prime$ of $A$ is $d+2$ and such a line bundle is necessarily globally generated. \\
Note that $d+2$ is the minimal degree of a line bundle with at least $3$ sections. If the vector bundle $E_{C,A}$ is simple, then we have $h^0(E_{C, A} \otimes {E_{C, A}}^*) =1 $. Thus 
 \begin{equation}\label{Q1}
  \chi(E_{C, A} \otimes {E_{C, A}}^*) = 2 -h^1(E_{C, A} \otimes {E_{C, A}}^*). 
 \end{equation}
 On the other hand, by Rieman-Roch, we have 
 \[
 \chi(E_{C, A}\otimes {E_{C, A}}^*) = \frac{{c_1(E_{C, A}\otimes {E_{C, A}}^*)}^2}{2} 
                            - c_2(E_{C, A} \otimes {E_{C, A}}^*) + \text{rk}(E_{C, A}\otimes {E_{C, A}}^*) \chi(\mathcal{O}_X).
 \]
Now $c_2(E_{C, A} \otimes {E_{C, A}}^*) = 6 c_2(E_{C, A}) - 2 {c_1(E_{C, A})}^2 $.
Thus we have, 
\begin{equation}\label{Q2}
\begin{split}
\chi(E_{C, A}\otimes {E_{C, A}}^*)\\
                  &= 18 - 6 c_2(E_{C, A}) +2 {c_1(E_{C, A})}^2 \\
                 & = 18-6(d+2) +2 (2g-2)\\
                 & =2-2\rho(g, 2, d+2).
\end{split}
\end{equation}                             
 Comparing \eqref{Q1} and \eqref{Q2} we have, $\rho(g, 2, d+2) \ge 0$. Thus $W^2_{d+2}(C)$ is non-empty for every smooth curve $C$ in $|L|$. Hence the second co-ordinate of the gonality sequence is constant.

 Let assume $E_{C, A}$ is not simple.   Then by Remark  \ref{RZ}, we have an exact sequence of the form 
 \begin{equation}\label{Q3}
 0 \to M \to E_{C, A} \to F \to 0,
 \end{equation}
 where $F$ is a rank $2$ torsion free sheaf, generated by its global sections off a finite set and  $M$ is line bundle, with at least two sections and $F$ fits in the following exact sequence,
 \[
 0 \to F \to F^{**} \to \mathcal{O}_Z \to  S \to 0,
 \]
 where $S$ is a coherent sheaf of finite length , in particular supported on a zero-dimensional subscheme $Z$.
 
  Let $N:=c_1(F)$.  
  Note   that 
 \begin{equation}\label{Q4}
c_2(E_{C, A})= d+2= M.N + |Z| + c_2(F^{**}).
 \end{equation}
 Since $F$ is globally generated by its section of a finite set, $F^{**}$  is also globally generated off a finite set.  Also note that as $F^{**}$ is globally generated by it's sections off a finite set, $N$ is globally generated off a finite set and since on a K3-surface a line bundle can have no fixed point outside  fixed components, $N$ is globally generated.\\ Also by Lemma \ref{L2}, $c_2(F^{**}) \ge 0$.\\
 Claim: $h^1(N) \le 1$.\\
 Since $N$ is base point free, $h^1(N) \ne 0$ implies that $N = \mathcal{O}(k\Gamma)$ [ Proposition 2.6 \cite{8}], where $\Gamma$ is an elliptic curve and $k$ is an integer $\ge 2$.   Also we have  $h^1(N) = k-1 \text{ and } h^0(N) = k+1 $. Thus if $h^1(N) > 1$, then $k \ge 3$.
  Since  $c_2(F) \ge 0$,  we have $C.2\Gamma < M.N \le d+2$ .  But $\mathcal{O}_C(2\Gamma)$ has 3 sections, which is a contradiction to the minimality of $d+2$.

In the case  when $h^1(N) =1$ we have, $N = \mathcal{O}(2\Gamma)$.  If $|Z| + c_2(F^{**}) > 0$, then $\text{deg}(N{\mid_C}) < d+2$  and $h^0(N{\mid_C})=3$. Thus we get a contradiction. \\
If $|Z| +c_2(F) = 0$, then $N{\mid_C}$ has $3$ sections  and degree of $N{\mid_C} = d+2$ for all $C \in |L|$, which proves our theorem.


Let us assume $h^1(N) = 0$.\\
 If $h^0(N) =2$, then $N = \mathcal{O}_X(E)$, where $E$ is an smooth elliptic curve.  On the other hand, since $F^{**}$ is globally generated off a finite set, by [\cite{4}, Proposition 1.5] , $F^{**} = \mathcal{O}_X(\Delta) \oplus \mathcal{O}_X(\Delta)$, where where $\Delta$ is a smooth
irreducible curve on $X$ which moves in a base-point free pencil. Thus $N = \mathcal{O}_X(2\Delta)$, a contradiction.

Let  $h^0(N) \ge 3$.  

 Since $h^0(M) \ge 2$ and $h^0(M) \le h^0(M{\mid_C}), M{\mid_C}$   contributes in the Clifford index. 
Since $K_C = \mathcal{O}_C(C)$, we have  $K_C \otimes {M_{\mid_C}}^* = N_{\mid_C}$.  From the exact sequence,
\[
0 \to \mathcal{O}(N-C) \to  N \to  N_{\mid_C} \to 0
\]
we have $h^0(N_{\mid_C}) = h^0(N) + h^1(M)$.  Also by Riemann-Roch, we have $h^0(N) = \frac{N^2}{2} + 2$. Thus 
 \begin{equation}\label{Q6}
 \begin{split}
\text{Cliff}(M{\mid_C})= \text{Cliff}(K_C \otimes {M{\mid_C}}^*)= \text{Cliff}(N{\mid_C})  &=  N.C -2 (h^0(N) + h^1(M)) +2 \\
                                      & = N.C - N^2 -4  -2h^1(M) +2\\
                                       &= M.N -2h^1(M) -2\\
                                        &=d+2 -|Z| -c_2(F^{**}) - 2h^1(M) -2.
\end{split}
\end{equation}
  But $d -2 =\text{Cliff}(C) \le \text{Cliff}(M{\mid_C})$, thus we have 

  \[
   d-2 \le d - |Z| -c_2(F^{**}) - 2h^1(M) 
   \]
   \[
\text{or }  |Z| + c_2(F^{**}) +2h^1(M) \le 2.
\]
 In particular $c_2(F^{**}) \le 2$.
 Since $F^{**}$ is globally generated off a finite set, for a general two dimensional subspace $V$ of $H^0(F^{**})$, we have 
\begin{equation}\label{Q5}
0 \to V\otimes \mathcal{O}_X \to F^{**} \to B \to 0
\end{equation}
where $B$ is a line bundle on a smooth curve $D \in |N|$.\\
Dualizing the above exact sequence we get,
\[
0 \to {F^{**}}^* \to V^* \otimes \mathcal{O}_X \to B^{\prime} \to 0
\]
where $B^{\prime} = \mathcal{O}_D(D) \otimes B^*$. Now from the long exact sequence of \eqref{Q3}, we have $h^0(F^*)=h^2(F) = 0$. Thus we have $h^0(B^{\prime}) \ge 2$. Also we have $c_2(F^{**}) = \text{deg}(B^{\prime})$.
  But $c_2(F^{**}) = \text{deg}(B^{\prime})$ and $B^{\prime}$ has at least 2 sections. Therefore the curve $D$ is hyperelliptic.  If $D$ has genus $2$ then $\text{deg}(\mathcal{O}(D){\mid_C})=D.C = D^2 + M.N= 2+M.N$. Since $c_2(F^{**}) =2, |Z| = 0$, then  from \ref{Q4}  it follows that $M.N = d$, which implies $D.C = d+2.$  Therefore $\mathcal{O}(D){\mid_C}$ will give a complete $g^2_{d+2}$ for all $C \in |L|$. If $D$ has genus bigger than $2$, then the following two cases can occur [\cite{8} Theorem 5.2] :\\
 (i) There exists an irreducible elliptic curve $\Delta$ such that $\Delta . D = 2$.\\
 (ii) There exists an irreducible hyperelliptic curve $B$ of genus $2$ such that  $D \sim 2B$.\\
 In case (i), we can further assume genus of $D$ is bigger than $3$,  thus we can   decompose $D$ as $\Delta + D^{\prime}$, with $D^{\prime}.\Delta = 2$. Now ${(D- 2 \Delta)}^2 = D^2 - 8$. Thus if $D-2\Delta$ is not effective, then $D^2 = 6$ and hence ${D^{\prime}}^2 = 2$.  Therefore the restriction of $\mathcal{O}(D^{\prime})$  on each curve in $|L|$ will give a complete $g^2_{d+2}$.  If $D-2\Delta$ is effective then  we can decompose $D$ as $D^{''} + 2 \Delta$ and $L = \mathcal{O}(2\Delta + D^{''}) \otimes M$. It is easy to see that ${(D^{''} + c_1(M))}^2 > 0$. Thus $D^{''}.c_1(M) \ge 2$, [\cite{8}, Lemma 3].\\ 
 On the other hand, 
 \begin{equation}
 \begin{split}
  \text{deg}(\mathcal{O}(2\Delta){\mid_C})=  4 + 2\Delta. c_1(M) \le M.N + c_2(F^{**}) +|Z| = d+2
 \end{split}
 \end{equation} 
 Therefore, $\mathcal{O}(2\Delta){\mid_C}$ will give a $g^2_{d+2}$ for all $C \in |L|$ or $\text{deg}(\mathcal{O}(2\Delta){\mid_C} < d+2$, a contradiction.
 
 In case (ii),  Considering the line bundle $\mathcal{O}(B)$. Note that since,  $C$ is neither  hyper-elliptic nor bi-elliptic, by Mumford's Theorem for $g^2_d$ \cite{1}, $W^2_{d+2}(C)$ is non-empty, if and only if, $d+2 - 6 \ge 0$ that is $d+2 \ge 6$, i.e., $d \ge 4$. \\
 Note that $M.N= d$  and 
  $B.C =B.(M+N)= B(2B +M)= 2B^2 + \frac{M.N}{2} \le d+ 2$. Thus either $\mathcal{O}_X(B){\mid_C}$ is a complete $g^2_d$ for all $C \in |L|$ or we will get a cotradiction.

 
{\bf Case II: r =3} 

Let $A$ be a line bundle of degree $d^{\prime}$  computing the Clifford index of $C$ with  $h^0(A) = 4$.
 We can assume there is no curve in $|L|$ with a line bundle with $3$ sections, computing the Clifford index.  Since $d^{\prime}$ computes the Clifford index of $C$ and the Clifford index of $C$ is $d-2$, one has  $d^{\prime} = d+4$. 
  In this case  every curve in the linear system $|L|$ admits a complete $g^3_{d+4}$ \cite[Theorem 4.1]{7}. 
For a general  point $x \in C$,  $A \otimes \mathcal{O}_C(-x)$ admits $3$ sections. Thus $W^2_{d+3}$ is non-empty. If $W^2_{d+2} \ne \o$, then one can get a line bundle computing the Clifford index of $C$ with $3$ sections, a contradiction. Thus  $W^2_{d+2} = \o$.
This is true for every smooth irreducible curve in $|L|$. Thus the planarity of every curve in the linear system $|L|$ is $d+3$.

{\bf Case III: r = 4}

Again let $A$ be a line bundle of degree $d^{\prime}$  computing the Clifford index of $C$ with  $h^0(A) = 5$. In this case $d^{\prime} = d + 6$ and as previous case by \cite[Theorem 4.1]{7}, every curve in the linear system $|L|$ admits a complete $g^4_{d+6}$. 
Now for   general  two points $x, y \in C, A \otimes \mathcal{O}_C(-x-y)$  admits $3$ sections. Thus $W^2_{d+4}(C)$ is non-empty for every smooth irreducible curve $C \in |L|$. If $W^2_{d+3}(C) = \o$ for all $C \in |L|$, then planarity of every curve is $d+4$ and we are done.

Let $C \in |L|$ such that $W^2_{d+3}(C) \ne \o$ and let $A \in W^2_{d+3}(C)$.\\
  Let $E_{C, A}, F^{**}, M, N , Z, D$ are as in Case I. Then from \ref{Q6}, we have 
\[
d-2 \le M.N -2h^1(M) -2 = d+3 - |Z| -c_2(F^{**}) -2h^1(M) -2
\]
Or
\[
c_2(F^{**}) + |Z| + 2h^1(M) \le 3
\]
If $c_2(F^{**}) \le 2$ then we can conclude the Theorem as Case I. Let $c_2(F^{**}) =3$. Then the degree of the line bundle $B$  on $D$ in \ref{Q5} is $3$ and admits $2$ sections. Hence $D$ is a trigonal curve. Since $L^2 \ge 16$ and $h^0(F) \ge h^0(N)$, we have $D^2 \ge 8$. In other words, $D$ has genus at least $5$. Therefore by Theorem \ref{TTA}, there exist an elliptic curve $\Delta$ is $X$ such that $\Delta . D = 3$ and $D$ can be decomposed as $D^{\prime} + \Delta$. Since $D^2 \ge 8$, we have ${D^{\prime}}^2 \ge 0$. If $0 \le {D^{\prime}}^2 \le 2$, then by similar analysis as in {\bf Case $r =2$}, we can conclude the Theorem. Thus we can assume that ${D^{\prime}}^2 \ge 4$, that is, $D^2 \ge 10$. If $D^2 \ge 12$, then $D$ can be decomposed as $2 \Delta + D^{\prime}$ and we are done as earlier. \\ 
Let $D^2 = 10$. Then ${D^{\prime}}^2 = 4$. Therefore, $D^{\prime}$ is either hyperelliptic or trigonal. Thus we have a decomposition of $D$ as $2\Delta + D^{\prime}$, which conclude the Theorem. 

  \end{proof}
  
  {\it Acknowledgement:} We would like to thank to Prof.A.J. Parameswaran for many useful discussion. We also would like to thank  Prof. Ciliberto , Prof. P. Newstead for  valueable 
comments and pointing out the work done in this direction. We also would like to thank Krishanu Dan 
for careful reading of the article.

 \end{document}